\documentclass[a4paper,11pt]{article}

\usepackage{amssymb,amsmath}
\usepackage{xypic}
\xyoption{all}
\usepackage{graphicx}

\newtheorem{defn}{Definition}[section]

\newtheorem{rem}[defn]{Remark}
\newtheorem{exm}[defn]{Example}

\newtheorem{notat}[defn]{Notation}
\newtheorem{newpar}[defn]{}

\newtheorem{xdefn}{Definition.}
\newtheorem{xproposition}{Proposition.}
\newtheorem{xcorollary}{Corollary.}
\newtheorem{xrem}{Remark.}
\newtheorem{xexm}{Example.}
\newtheorem{xlemma}{Lemma.}
\newtheorem{xtheorem}{Theorem.}
\newtheorem{xnotat}{Notation.}
\newtheorem{xnewpar}{\it}
\newtheorem{xproof}{{\it Proof. }}
\newtheorem{xproofof}{{\it Proof}}

\newenvironment{proof}{\begin{xproof}\em}{\end{xproof}}

\newenvironment{newparagraph*}[1]{\begin{xnewpar}\hspace*{-1.5mm}{#1}. \rm}{\end{xnewpar}}

\newenvironment{definition*}{\begin{xdefn}\em}{\end{xdefn}}
\newenvironment{remark*}{\begin{xrem}\em}{\end{xrem}}
\newenvironment{example*}{\begin{xexm}\em}{\end{xexm}}
\newenvironment{notation*}{\begin{xnotat}\em}{\end{xnotat}}
\newenvironment{proposition*}{\begin{xproposition}}{\end{xproposition}}
\newenvironment{corollary*}{\begin{xcorollary}}{\end{xcorollary}}
\newenvironment{lemma*}{\begin{xlemma}}{\end{xlemma}}
\newenvironment{theorem*}{\begin{xtheorem}}{\end{xtheorem}}

\def\qed{\hspace{0.3cm}{\rule{1ex}{2ex}}}
\newcommand\V{\bigvee}

\newcommand\ie{i.e.}

\newcommand\st{\mid}

\newcommand\downsegment{{\downarrow}}

\begin{document}

\title{A note on infinitely distributive inverse semigroups%
\thanks{Research supported in part by FEDER
and FCT through CAMGSD.}}
\author{Pedro Resende
\vspace*{2mm}\\ \small\it Departamento de Matem{\'a}tica, Instituto
Superior T{\'e}cnico, \vspace*{-1mm}\\ \small\it Av. Rovisco Pais 1,
1049-001 Lisboa, Portugal}

\date{~}

\maketitle

\noindent By an \emph{infinitely distributive} inverse semigroup will be meant an inverse semigroup $S$ such that for every subset $X\subseteq S$ and every $s\in S$, if $\V X$ exists then so does $\V (sX)$, and furthermore $\V(sX)=s\V X$.

One important aspect is that the infinite distributivity of $E(S)$ implies that of $S$; that is, if the multiplication of $E(S)$ distributes over all the joins that exist in $E(S)$ then $S$ is infinitely distributive. This can be seen in Proposition 20, page 28, of Lawson's book \cite{Lawson}. Although the statement of the proposition mentions only joins of nonempty sets, the proof applies equally to any subset.

The aim of this note is to present a proof of an analogous property for binary meets instead of multiplication; that is, we show that for any infinitely distributive inverse semigroup the existing binary meets distribute over all the joins that exist.

A useful consequence of this lies in the possibility of constructing, from infinitely distributive inverse semigroups, certain quantales that are also locales (due to the stability of the existing joins both with respect to the multiplication and the binary meets), yielding a direct connection to \'{e}tale groupoids via the results of \cite{Resende}. The consequences of this include an algebraic construction of ``groupoids of germs'' from certain inverse semigroups, such as pseudogroups, and will be developed elsewhere.

\begin{lemma*}
Let $S$ be an inverse semigroup, and let $x,y\in S$ be such that the meet $x\wedge y$ exists. Then the join
\[f=\V\{g\in \downsegment (xx^{-1}yy^{-1})\st gx=gy\}\]
exists, and
we have $x\wedge y=fx=fy$.
\end{lemma*}

\begin{proof}
Consider the set $Z$ of lower bounds of $x$ and $y$,
\[Z=\{z\in S\st z\le x,\ z\le y\}\;,\]
whose join is $x\wedge y$. By \cite[Prop.\ 17, p.\ 27]{Lawson}, the (nonempty) set
\[F=\{zz^{-1}\st z\in Z\}\]
has a join $f=\V F$ that coincides with $(\V Z)(\V Z)^{-1}=(x\wedge y)(x\wedge y)^{-1}$. Hence, $x\wedge y=fx=fy$. The lemma now follows from the fact that the elements $zz^{-1}$ with $z\in Z$ are precisely the idempotents $g\in\downsegment(xx^{-1}yy^{-1})$ such that $gx=gy$.
\qed
\end{proof}

Under the assumption of infinite distributivity we have a converse:

\begin{lemma*}
Let $S$ be an infinitely distributive inverse semigroup, and let $x,y\in S$ be such that the join
\[f=\V\{g\in \downsegment (xx^{-1}yy^{-1})\st gx=gy\}\]
exists. Then the meet $x\wedge y$ exists, and we have $x\wedge y=fx=fy$.
\end{lemma*}

\begin{proof}
By \cite[Prop.\ 20, p.\ 28]{Lawson}, the join
\[\V\{gx\st g\le xx^{-1}yy^{-1},\ gx=gy\}\]
exists and it equals $fx$. Similarly, the join
\[\V\{gy\st g\le xx^{-1}yy^{-1},\ gx=gy\}\]
exists and it equals $fy$. But the two sets of which we are taking joins are the same due to the condition $gx=gy$, and thus $fx=fy$. The element $fx$ is therefore a lower bound of both $x$ and $y$. Let $z$ be another such lower bound. Then $z=zz^{-1}x=zz^{-1}y$, and thus $zz^{-1}\le f$, which implies $z\le fx$. Hence, $fx$ is the greatest lower bound of $x$ and $y$.
\qed
\end{proof}

\begin{theorem*}
Let $S$ be an infinitely distributive inverse semigroup, let $x\in S$, and let $(y_i)$ be a family of elements of $S$. Assume that the join $\V_i y_i$ exists, and that the meet $x\wedge\V_i y_i$ exists. Then,
for all $i$ the meet $x\wedge y_i$ exists, the join $\V_i(x\wedge y_i)$ exists, and we have
\[x\wedge\V_i y_i = \V_i(x\wedge y_i)\;.\]
\end{theorem*}

\begin{proof}
Let us write $y$ for $\V_i y_i$, $e_i$ for $y_i y_i^{-1}$, and let $f$ be the idempotent
\[f=\V\{g\in \downsegment(xx^{-1}yy^{-1})\st gx=gy\}\;,\]
which exists, by the first lemma. Furthermore, also by the first lemma, we have $x\wedge y=fx=fy$.
We shall prove that $x\wedge y_i$ exists for each $i$, and that it equals $e_i (x\wedge y)=e_i fx=e_i fy$.
By the second lemma, it suffices to show that for each $i$ the join
\[f_i=\V\{g\in \downsegment(xx^{-1}y_i y_i^{-1})\st gx=gy_i\}\]
exists and equals $e_i f$. Consider $g\in \downsegment(xx^{-1}y_i y_i^{-1})=\downsegment(xx^{-1} e_i)$ such that $gx=gy_i$.
The condition $g\in\downsegment(xx^{-1}e_i)$ implies that $g\le e_i$, and thus $g=g e_i$. Hence, since
$y_i=e_i y$, the condition $gx=gy_i$ implies $gx=ge_i y=gy$, and thus $g\le f$ because furthermore
$g\in\downsegment(xx^{-1}yy^{-1})$. Hence, we have both
$g\le e_i$ and $g\le f$, \ie, $g\le e_i f$, meaning that $e_i f$ is an upper bound of the set
\[X=\{g\in \downsegment(xx^{-1}e_i)\st gx=gy_i\}\;.\]
In order to see that it is the least upper bound it suffices to check that $e_i f$ belongs to $X$, which is immediate: first, $e_i f\le xx^{-1}$ because $f\le xx^{-1}$, and thus $e_i f\in\downsegment( xx^{-1} e_i)$; secondly,
\[(e_i f) x= (e_i  e_i f) y = (e_i f e_i) y = (e_i f)(e_i y) = (e_i f) y_i\;.\]
Hence, $e_if\in X$, and thus $e_i fx=x\wedge y_i$.
In addition, the join $\V_i e_i$ exists and it equals $yy^{-1}$, by \cite[Prop.\ 17, p.\ 27]{Lawson}, and thus, using infinite distributivity and the fact that $f\le yy^{-1}$, we obtain
\[x\wedge y=fx=yy^{-1} fx=(\V_i e_i) fx=\V_i (e_i f x)=\V_i(x\wedge y_i)\;. \qed\]
\end{proof}

\end{document}